\documentclass[11pt, a4paper]{amsart}

\usepackage{amsmath}
\usepackage{amstext}
\usepackage{amssymb}
\usepackage{amsthm}
\usepackage{enumerate}
\usepackage[shortlabels]{enumitem}
\usepackage{mathrsfs}
\usepackage[all]{xy}
\usepackage{datetime}

\usepackage[svgnames]{xcolor}
\usepackage{tikz}

\pgfdeclarelayer{edgelayer}
\pgfdeclarelayer{nodelayer}
\pgfsetlayers{edgelayer,nodelayer,main}

\tikzstyle{none}=[inner sep=0pt]

\makeatletter
\def\imod#1{\allowbreak\mkern10mu({\operator@font mod}\,\,#1)}
\makeatother

\newtheorem{theorem}{Theorem}[section]
\newtheorem{prop}[theorem]{Proposition}
\newtheorem{lemma}[theorem]{Lemma}

\newtheorem{corollary}[theorem]{Corollary}

\theoremstyle{definition}
\newtheorem{definition}[theorem]{Definition}

\newtheorem{question}{Question}

\theoremstyle{remark}
\newtheorem{remark}[theorem]{Remark}

\theoremstyle{remark}

\numberwithin{equation}{section}

    \DeclareMathOperator{\lh}{lh}

    \DeclareMathOperator{\ran}{ran}

    \DeclareMathOperator{\fin}{FIN}
    
    \DeclareMathOperator{\card}{card}

    \DeclareMathOperator{\hod}{HOD}
    \DeclareMathOperator{\od}{OD}
    \DeclareMathOperator{\add}{add}
        \DeclareMathOperator{\bp}{BP}
        \DeclareMathOperator{\lm}{LM}
        \DeclareMathOperator{\nd}{NDe}
        \DeclareMathOperator{\cd}{CDi}                
        \DeclareMathOperator{\interior}{Int}                
    
    \newcommand{\restrict}{\!\upharpoonright\!}

    \newcommand{\append}{^\frown}

\newcommand{\uT}{\mathfrak u_T}
\newcommand{\uM}{\mathfrak u_M}
\newcommand{\calT}{\mathcal T}
\newcommand{\emptyinfinite}{\{0,\aleph_0\}}

\def\N{{\mathbb N}}

\def\P{{\mathbb P}}

\def\Q{{\mathbb Q}}

\title[Set theory and a model of the mind in psychology]{Set theory and a model of the mind in psychology}

\author{Jens Mammen, Asger T\"ornquist}

\address{Department of Communication and Psychology, Aalborg University, Rendsburggade 14, 9000 Aalborg, Denmark. Orcid: 0000-0002-8867-5744.}

\email{mammen@hum.aau.dk}

\address{Department of Mathematical Sciences, University of Copenhagen, Universitetsparken 5, 2100 Copenhagen, Denmark. Orcid: 0000-0002-8332-649X.}

\email{asgert@math.ku.dk}


\date{\today}

\begin{document}

\maketitle

\begin{abstract}
We investigate the mathematics of a model of the human mind which has been proposed by the psychologist Jens Mammen. Mathematical realizations of this model consist of so-called \emph{Mammen spaces}, where a Mammen space is a triple $(U,\mathcal S,\mathcal C)$, where $U$ is a non-empty set (``the universe''), $\mathcal S$ is a perfect Hausdorff topology on $U$, and $\mathcal C\subseteq\mathcal P(U)$ together with $\mathcal S$ satisfy certain axioms.

We refute a conjecture put forward by J. Hoffmann-J{\o}rgensen, who conjectured that the existence of a ``complete'' Mammen space implies the Axiom of Choice, by showing that in the first Cohen model, in which ZF holds but AC fails, there is a complete Mammen space. We obtain this by proving that in the first Cohen model, every perfect topology can be extended to a maximal perfect topology.

On the other hand, we also show that if all sets are Lebesgue measurable, or all sets are Baire measurable, then there are no complete Mammen spaces with a countable universe.

Finally, we investigate two new cardinal invariants $\uM$ and $\uT$ associated with complete Mammen spaces and maximal perfect topologies, and establish some basic inequalities that are provable in ZFC. Further, we show $\uM=\uT=2^{\aleph_0}$ follows from Martin's Axiom, and, contrastingly, we show that $\aleph_1=\uM=\uT<2^{\aleph_0}=\aleph_2$ in the Baumgartner-Laver model.
\end{abstract}

\section{Introduction}

In theoretical psychology, Jens Mammen has proposed a model for what may be called the \emph{interface} between the inner world of a human mind, and the outer world that this human mind lives in, perceives, and interacts with. From the outset, Mammen has formulated and presented his theory axiomatically, in the style familiar to mathematicians. The purpose of this paper is to study the set-theoretic aspects of Mammen's theory.

\medskip

Briefly, a Mammen space can be defined as follows:

\begin{definition}\label{d.mammen} A Mammen space is a triple $(U,\mathcal S,\mathcal C)$, where $U$ is a non-empty set, called the \emph{universe of objects}, and $\mathcal S,\mathcal C\subseteq \mathcal P(U)$ such that
\begin{enumerate}
\item $\mathcal S$ is a perfect Hausdorff topology on $U$; here \emph{perfect} means that every non-empty open set is infinite, and so in particular the set $U$ is infinite.
\item $\mathcal C$ satisfies:
\begin{enumerate}
\item There is a non-empty $C\in\mathcal C$;
\item $\mathcal C$ is closed under finite unions and intersections;
\item Every non-empty $C\in\mathcal C$ contains a singleton which is in $\mathcal C$.
\end{enumerate}
\item $\mathcal S$ and $\mathcal C$ together must satisfy
\begin{enumerate}
\item $\mathcal S\cap \mathcal C=\{\emptyset\}$;
\item If $C\in\mathcal C$ and $S\in\mathcal S$ then $C\cap S\in\mathcal C$.
\end{enumerate}
\end{enumerate}
\end{definition}
The elements of $\mathcal S$ are called \emph{sense categories} and the elements of $\mathcal C$ are called \emph{choice categories}. The reader should think of a Mammen space $(U,\mathcal S,\mathcal C)$ as a model that a person's mind has (or has built) as a result of sensory input and experience: It has formed broad categories of the objects in the universe $U$, and these are represented by the subsets of $U$ which are in $\mathcal S$; and it has singled out categories of particular objects or people, and these are represented by subsets of $U$ which are in $\mathcal C$. For instance, the mind of a person overlooking a beach will have a sense category of all the stones on the beach, but if that person singles out a special stone and picks it up, he is availing himself of a choice category, which in the mathematical representation is the singleton of that special stone. Restating this with the emphasis on the role of $\mathcal C$ instead, the idea is that categories in $\mathcal C$ represent collections of objects, people, animals, etc., of particular attachment for the person (e.g., the person's father), in contrast to the broad categories in $\mathcal S$ (e.g., the category of all people who are fathers).

Fuller details of the psychological background and motivation for the definition of a Mammen space is given in section \ref{s.psych} below.

\medskip

The question which gives rise to much of the mathematics of this paper is the question of \emph{completeness}: Are the categories $\mathcal S$ and $\mathcal C$ sufficient to be able to account for \emph{all} possible categories of objects that can be formed in the universe? That is, can every $X\subseteq U$ can be written as
$$
X=S\cup C\text{ where } S\in\mathcal S\text{ and } C\in\mathcal C\text{\ \ ?}
$$
If this is the case, we will call the Mammen space $(U,\mathcal S,\mathcal C)$ \emph{complete}. 

The question of existence of a complete Mammen space turns out to be mathematically non-trivial. It was answered in the positive by J. Hoffmann-J{\o}rgensen in \cite{HJ2000}, but the Axiom of Choice (below abbreviated AC, or Choice) was used to do so:

\begin{theorem}[Hoffmann-J{\o}rgensen. Uses AC]\label{t.hoffmannintro}
For any infinite set $U$, there is a complete Mammen space with universe $U$.
\end{theorem}

Hoffmann-J{\o}rgensen proved this by observing that if $\mathcal S$ is a \emph{maximal perfect topology} on $U$, and we take $\mathcal C$ to be the family of closed nowhere dense subsets of $U$, then $(U,\mathcal S,\mathcal C)$ is a complete Mammen space. (We will reprove this below, see Theorem \ref{t.hoffmann} and Corollary \ref{c.hoffmann}.) The Axiom of Choice is used by Hoffmann-J{\o}rgensen only to ensure that maximal perfect topologies exist.

\medskip

Do we really need the Axiom of Choice to prove Theorem \ref{t.hoffmannintro}? Given the psychological origin and relevance of the notion of a complete Mammen space, it is desirable to avoid using AC in its full strength, if possible, and at the same time highly interesting if AC can't be avoided. Hoffmann-J{\o}rgensen was of the opinion that AC is unavoidable, and conjectured this:

\medskip

\noindent {\bf Conjecture I.} (Hoffmann-J{\o}rgensen).  
If there is a complete Mammen space, then the Axiom of Choice holds. In particular, if there is a maximal perfect Hausdorff topology, then the Axiom of Choice holds.
\medskip

The purpose of this paper is to refute this conjecture, and at the same time show that Theorem \ref{t.hoffmannintro} requires some non-constructive mathematical methods, at least if the universe $U$ is required to be countably infinite. Specifically, we will prove:

\medskip

\noindent {\bf Theorem A.}\label{t.mainintro1}\emph{
In the first Cohen model, that is, $\hod^{V[G]}(A)$ where A is the countable set of Cohen reals added by the generic $G$, there is a complete Mammen space, and the underlying universe $U$ of this space can even be chosen to be countable. Indeed, in $\hod^{V[G]}(A)$, it holds that every perfect topology can be extended to a maximal perfect topology.}

\medskip

Since it is well-known that $\hod^{V[G]}(A)$ is a model of Zermelo-Fraenkel (ZF) set theory in which Choice \emph{fails}, the previous theorem refutes Hoffmann-J{\o}rgensen's conjecture.

\medskip

As a counterpoint to Theorem A, we will show:

\medskip

\noindent{\bf Theorem B.} \emph{If all sets are Lebesgue measurable, or if all sets are Baire measurable, then there are no complete Mammen spaces with a countable universe $U$. It follows that it is not possible to prove in ZF alone that there is a complete Mammen space with a countable universe.}

\medskip

In the final parts of the paper, we will also examine the following two cardinal invariants that are naturally associated with the notions of maximal perfect topologies and complete Mammen spaces:
$$
\uT=\inf\{\card(\mathcal T): \mathcal T \text{ is a maximal perfect Hausdorff topology on }\N\}.
$$
and
$$
\uM=\inf\{\card(\mathcal S): \mathcal S \text{ is a sense category of a complete Mammen space on }\N \}.
$$
We will show the following:

\medskip

\noindent {\bf Theorem C.} {\it Let $\add(\bp)$ denote the additivity of the meagre ideal, see \cite[p. 515]{Jech03}. Then:

\begin{enumerate}
\item $\add(\bp)\leq\uM\leq \uT\leq 2^{\aleph_0}$
\item In particular, $\uM$ and $\uT$ are always uncountable, and if Martin's Axiom holds then $\uM=\uT=2^{\aleph_0}$.
\item In the Baumgartner-Laver model we have $\aleph_1=\uM=\uT<2^{\aleph_0}=\aleph_2$.
\end{enumerate}
}

Towards the end of the paper, we discuss several intriguing questions that remain open.

\section{Psychological background}\label{s.psych}

The present paper is basically about some mathematical problems. However, the motivation for posing the problems has its roots in corresponding questions in psychological science. Therefore, there will be a short introduction to these psychological questions, also to provide a possible interpretative frame for the mathematics, or even a model outside mathematics. At the same time, we attempt to briefly show the relevance to psychology of the kind of mathematical model that is presented here.

This introduction is not in itself a psychological paper with the usual demands for documentation and direct references, as this would consume too much space in the present context. More details could be found in the works in the reference list.

\subsection{Introducing psychology and psychophysics}

Today, psychology is not a coherent science with commonly accepted basic theoretical concepts. This means that it will be impossible to give a short covering definition of the scientific field, and perhaps even worse, also of the concrete domain of study, or in other words what it is about and how it can be applied as a tool in this context.

There is, however, a not quite negligible minority claiming that human psychology must in some way be about the ``interface''\footnote{The term ''interface'' will be used for the ''practical'' or ``active'' relation between humans and the world of objects. It is not just referring to a surface of contact but to a relational structure expanding in space and time. Perhaps ``interspace'' would have been more precise, as proposed by Engelsted \cite{Engelsted}. Here we have, however, chosen to follow conventional terminology.} between humans and the world we are living in.

This does not mean that questions of what is going on inside the body, and especially in the brain, be it subjective experiences and/or physiological processes, is of no interest in psychology, on the contrary. But a prerequisite for this study is that the tasks to be solved by the brain and the body meeting the world are rather well understood. Walking is the key to understanding the legs, which again serve and constrain walking.

The mythical, philosophical, and scientific understanding of this ``interface'' has a long history since antiquity which of course can't be covered here. When focusing on scientific psychology something dramatic happened, however, around 1850 when psychology found a way to define itself as a natural science in the conceptual frame of contemporary physics, chemistry, physiology and mathematics. Before that psychology had rather been considered an auxiliary discipline to theology and philosophy.

It is thus common to define the birth of modern scientific psychology to the introduction of so-called psychophysics, often referring to the theory of G. T. Fechner. The idea was to consider the senses, e.g. vision, hearing, smelling, as ``transmitters'' receiving objective physical-chemical ``input'', and as ``output'' causing subjective impressions with some, mostly hypothetical, physiological  correlates or equivalents. Further the idea was to apply measures not only to the objective input but also to the subjective output, so the two events or entities could be bridged by a quantitative mathematical function.

This bridging, however, had a price. When our contact with the world before psychophysics had been understood as meaningful it was now reduced to the raw material of patterns of, in themselves meaningless, quantitative sense impressions. How could meaning, or conceptual knowledge, be reestablished on this meager ground? Psychology was ``eaten'' by the mechanistic understanding of the outer world. Psychophysics not only appeared as a bridge but perhaps even more as a barrier between man and world.

The problem is classic and reflected in European philosophy since the renaissance. There are inductive or empiricist attempts appealing to high, complex or hierarchical organizations of input-patterns (``sense data'') hoping that ``consciousness'' would pop up with enough complexity, in vain of course. And there are deductive or rationalist attempts appealing to {\it a priori} conceptual frames inducing order and meaning in the patterns, e.g. as when I. Kant rightly claims that time and space as frames for objects can't be inferred from sense impressions but have to be \emph{a priori}, but just raising new problems.

Many other attempts have been promoted to overcome the reductions of psychophysics. There has been appeal to language, hermeneutics, and semiotics in what in newer philosophy and psychology has been called ``the linguistic turn''. There are even attempts to reintroduce Aristotelian teleology violating the modern concept of proximal forward causality.

The result is that today we either have a reductionist mechanistic psychology or a psychology with a schism between a pure naturalist approach and a pure humanistic approach, sometimes expressed as ``the two cultures'' or \emph{Naturwissenschaft} versus \emph{Geisteswissenschaft} with two incompatible frames of understanding, causing both theoretical and practical problems.

Common to these attempts or approaches is that they don't correct or change psychophysics but either accept it as it is, or try to supplement it with principles ``taken from elsewhere'' but basically being incompatible with psychophysics. A third stance is just to turn your back to psychophysics and natural science and promote a pure humanistic psychology.

\subsection{Inspiration from modern natural science and mathematics}

There is, however, still another approach with inspiration from physics. When it was discovered that electromagnetic propagation of waves and particles did not follow the same kinematic laws as movements of solid bodies you did not choose a split in theory of movements in space and time but rather, as Einstein did, searched a ``conservative generalization'' of classical kinematics to include both corporeal movements and electromagnetic movements in one common law.

Einstein analyzed in detail the classic ``Galilei-transformation'' for movements of bodies in different systems of reference and searched what was the minimal change which conserved the classical laws for slow movements and still included the new knowledge of the speed of light. The apparent paradox is that this conservative approach implied the most revolutionary result, the famous formula for equivalence of mass and energy, as a simple deductive implication.

The principle of conservative generalization is also well known from mathematics:  Non-Euclidean geometries are still locally Euclidean, complex numbers include the real numbers, etc. The principle has, however, not played an apparent role in psychology.

Let us have a closer look at psychophysics to see if we can copy Einstein and conserve its unquestionable virtues without falling back to its reductionism. First, we must re-conceptualize psychophysics.

\medskip

Early psychophysicists believed that they bridged objective stimulation and subjective impressions when measuring both and connecting them in a quantitative functional relationship.

That was, however, result of a rather speculative interpretation of what was going on in the experiments. What was demanded from the experimental subjects when they were presented for physical stimulation was in fact only to make a yes/no-decision, or with other words to react or not. The questions to answer were either if some stimulation was at all being noticed (being greater than zero), or if a stimulation in some well-defined respect (e.g. size, strength, pitch) was greater than (or smaller than) another stimulation, serving as a comparative standard.

The experiments can be considered a continuous mapping from the domain of stimulation on the two-valued set (yes, no), and what was found in the experiments was the inverse images in the domain of the response ``yes'' (for noticed difference).\footnote{If the domain of stimulation is a continuum there will, however, necessarily be some uncertainty in the responses, so a more realistic model of the situation will be to allocate a probability for responding ``yes'' to each point in the stimulus domain.
In this case, what is found in the experiments is rather the inverse images in the domain of stimulation of the mapping on open intervals in the set of probabilities as e.g. ``points in the stimulation set with a probability greater than 0.5 for the yes-response''.
This is a generalization \cite{Dzhafarov} of the non-probabilistic case described above, which accordingly can be seen as a limiting case with probabilities 0 for the comparative standards themselves and 1 for all other points.
As the mathematical consequences in the present context are the same as for the non-probabilistic case we shall not dwell more on this generalization, well-known to psychologists.
}

These kinds of inverse images were introduced in \cite[1983]{Mammen96}, and following the terminology established there, they are called \emph{sense categories}.  They are sets organizing the domain of stimulation in a structure similar to the way the real axis is organized by open intervals defined by measurements based on the relations ``greater than'' and ``smaller than'', which is a perfect Hausdorff topology or just a perfect topology\footnote{This not only holds true for ``the naked senses'', as traditionally studied in psychophysics, but also for the senses ``expanded'' with amplifying tools and measuring equipment.}.

It is postulated that also outside the experimental situation is this topology a structural description of the senses' capacity, of course to be filled out with more quantitative definitions. Still this may be an idealization, but perhaps the best one we have as a theoretical foundation for understanding our sensory interface with objects in our world.

This description of our ``interface'' with domains of stimulation from objects in our environment can be generalized to the domain of the objects themselves. In psychology this is conceptualized as the movement, or step, from sensation to perception.

Now, the sense categories can be seen as organizing the domain of objects in a structure equivalent to a perfect topology.

Furthermore, this step from sensation to perception includes perception in a more general approach in psychology which could be called the extensional approach, i.e. the attempt to understand the human subject by taking departure in which parts of the, in principle infinite, world we are making objects for our relations and activities. In this context ``objects'' should be understood as including places as well as other subjects.

A little metaphorically expressed are our practical, cognitive and emotional relations and activities considered as selections or ``figure-ground'' operations on the infinite domain of objects, initiated by the subject. The originally Russian Activity Theory \cite{Leontiev} is a paradigmatic example of thinking about subjects in terms of such object directed activities. The so-called ecological approach as represented by the American psychologist J. J. Gibson \cite{Gibson} is another example.

\subsection{A conservative generalization of psychophysics}

In this perspective it becomes evident what is the insufficiency of reducing humans' ``interface'' with the world to a structure of sense categories. Sense categories are general categories of measures or ``universals'' and are only catching in principle infinite sets of objects defined by their measurable properties. They are not able to ``zoom in'' on any particular object. The perfect topology has no singletons.

But humans are not only living in a world of ``superficial'' universal properties, corresponding to the perspective of Artificial Intelligence. We are first of all relating to ``particulars'' or ``individuals'' with an individual history, in many cases irreplaceable and linguistically denoted by proper names. That is the case with our relatives, our possessions, and our belonging. It is these ``deep'' relations of co-existence which give our life meaning, and it is fatal if psychology, of all sciences, is ignoring them.

The historical ``threads'' of particular objects are also what define the meanings we share and express in language and concepts as e.g. a present or a gift. The difference between a valid coin and a counterfeit is not their properties but their individual history of origin.

But the ignorance is also fatal in a practical sense. As already Kant pointed out is the condition for an empirical statement (not only in science but also in everyday life) that the chosen particular object of predication is defined independently of the universal predicates in the statement. If you already have used them for identification, the statement is not empirical (synthetic) but analytic.

Our \emph{choice} of objects, in space and time, for use or investigation, or predication, is not dependent of an infinite process of ``zooming in'' on the set of universal sense categories. Due to \emph{our existence as particulars ourselves} and being in a particular place at a particular time we can just take an object or point it out.

\medskip

When walking on the beach I can just pick up an accidental stone without having to define it in advance by discriminating it from all other stones. I can put it in my pocket, and without having noticed its form and color I can be sure it is the same when returning home.

\medskip

In contrast to sense categories, particulars or collections of particulars are here called \emph{choice categories} (following the terminology established in \cite[1983]{Mammen96}). They are not necessarily finite, as they could also be defined by networks departing from particulars as a, in principle infinite, genealogical tree.

\medskip

These two structures are disjoint in the sense that no non-empty choice category can be a sense category, although they of course may share objects. But at the same time the two structures are framing each other. When picking up a stone I am not searching a piece of driftwood. And when coming home with a finite collection of chosen stones I will be able to not only distinguish them mutually but also to identify each -- within this collection -- with a finite sensory description.

\medskip

In fact, this capacity to have simultaneous dual relations to objects in the world, as members of sense categories and of choice categories, may have some antecedents in animals, but seems, in its full realization, to be a human privilege. In philosophical terms it is the capacity to operate jointly with objects' qualitative as well as numerical identity. Besides being basis for establishing a genuine referential language, it makes the distinction meaningful between e.g. seeing a new object, and seeing a well-known object with changed position or properties, which is vital for our cognitive, practical and emotional life.

This capacity is developed during our first year of life, before our appropriation of language, and remains a logical basis also in adult life. A renowned experimental study in this context is \cite{Xu}. For an overview of some later research see \cite{kroejgaard}.

\subsection{Axiomatics for sense and choice categories}

It should now be time for presenting \textit{an axiomatic system} describing the \textit{joint} structure of sense and choice categories, expressed in first-order equivalent language. Here $U$ denotes the world of \textit{objects}. 

\textbf{Ax. 1}: There is more than one object in $U$

\textbf{Ax. 2}: The intersection of two sense categories is a sense category 

\textbf{Ax. 3}: The union of any set of sense categories is a sense category 

\textbf{Ax. 4} (Hausdorff): For any two objects in $U$ there are two disjunct sense categories so that one object is in the one and the other object in the other one 

\textbf{Ax. 5} (perfectness): No sense category contains just one object 

\textbf{Ax. 6}: No non-empty choice category is a sense category 

\textbf{Ax. 7}: There exists a non-empty choice category 

\textbf{Ax. 8}: Any non-empty choice category contains a choice category containing only one object 

\textbf{Ax. 9}: The intersection of two choice categories is a choice category 

\textbf{Ax. 10}: The union of two choice categories is a choice category 

\textbf{Ax. 11}: The intersection of a choice category and a sense category is a choice category

It has been proven that the axioms are consistent and independent \cite{Mammen96,Mammen17}.

Axioms Ax. 1--5 state that sense categories are the open sets in a perfect topology on the underlying set of objects $U$, and so correspond to (1) of Definition \ref{d.mammen} in the introduction.

Axiom Ax. 5 claims that there are no singletons or that no single object is ``decidable'' in the topology of sense categories.

If $U$ had been finite would axiom Ax. 4 imply that all single objects themselves were sense categories, or singletons, in the topology. This is however negated by Ax. 5, which proves, that $U$ must be infinite.

It is Ax. 5 which ``opens for'' or ``makes room for'' the existence of non-empty choice categories as stated in Ax. 7 and thus invites the conservative generalization of the topology of sense categories to a structure also including choice categories.\footnote{This generalization also means that sense categories are only bound by the axioms and no longer by some order rooted in the ``greater than'' or ``smaller than'' relations in the special case of psychophysics, chosen here as a historical point of departure, but also to connect with a well-known field of psychology already formalized mathematically. In \cite{Mammen17} the same axioms are introduced more directly and logically from a concept of decidability and independent of psychophysics. Finally, in the 1983 version of \cite{Mammen96}, the 11 axioms are introduced as generalizations of the case with a finite $U$, corresponding to the typical case in experimental cognitive psychology. 
It may be a point of its own that the same axiomatic structure can be reached in at least three different ways.
}

Axiom Ax. 6 states the mutual exclusion of the two kinds of categories, and correspond to (3.a) of Definition \ref{d.mammen} in the introduction.

Axioms Ax. 7-10 describe the structure of choice categories, and correspond to part (2) of Definition \ref{d.mammen}.

Axiom  Ax. 8 secures the existence of finite non-empty choice categories. Further it states, in interpretative terms, that every non-empty choice category contains an ``accessible'', ``reachable'' or ``decidable'' member, or in other words, that every non-empty choice category must contain at least one identified specimen or instance. It further follows from the axioms that after picking out a member of a choice category, what is left is still a choice category. But also that it does not follow that the choice category necessarily could be ``emptied'' or ``exhausted'' by repeating this operation. It does also not follow that every member of a choice category is a choice category. That would be too radical generalizations, although the axioms don't exclude these possibilities.

Finally, is axiom Ax. 11 expresses the interaction, or mutual framing, of the two kinds of categories. This corresponds to (3.b) of Definition \ref{d.mammen}.

Of course, we can also combine the categories and define a joint concept of decidable category:

\begin{definition} 
A decidable category is a union of a sense category and a choice category.
\end{definition}

As it can be proven from the axioms that the empty set $\emptyset$ is both a sense category and a choice category, it follows that sense and choice categories themselves are decidable categories.

There are many interpretations and implications of the axiomatic system for psychological theories and their application, e.g. clinical psychology and developmental psychology. But here we have to refer to the reference list.\footnote{Still, the interpretation presented until now is, in relation to psychology, an idealization. Instead of considering the domain of stimulations as static ``properties'' being ``measured'' by comparison with objective or subjective ``standards'', it is rather variables being continuous functions of time and/or of explorative actions. An explorative action could e.g. be a pressure applied on some object which together with the resulting deformation as a function of the variable pressure would provide sensory information about the object's elasticity. Or it could be the way we often actively ``rock'' an object in our hand to provide not only sensory information of its mass as gravitational weight but also as inertial resistance to variable acceleration. The sense categories now become inverse images of responses in this generalized domain of functions.

Correspondingly, the objects in $U$ themselves, as also being defined by their individual history, could rather be interpreted as continuous ``threads'' in space and time. 

These expansions are parallel to the generalization of a point-topology to a topology of continuous functions as in the compact-open topology for function spaces \cite{Fox}; \cite[pp. 83-84]{Mammen17}). The explorative actions could further be expanded with explorations mediated by tools as microscopes and even chemical analyses.
The structure of sense and choice categories described in the 11 axioms in fact seems to be invariant to all these expansions (\cite[pp. 83--84]{Mammen17}). In other words, it is hypothesized that the dual structure of sense and choice categories is pervasive in a broad field of interpretative expansions, or applications, which were omitted in the introduction for simplicity reasons, and also to be closer to the historical development of modern psychology from psychophysics.

As the focus in the present paper is on this basic structure in itself, and not primarily on the interpretative fields, we have no intentions to cover those, except with a few chosen examples.
}

Here shall just be referred to two consequences in form of derived theorems, Th. 9 and Th. 10, from the set of axioms. The numbering refers to the one used in \cite[pp. 80--82]{Mammen17}.

\begin{theorem}[Th. 9, Correspondence] Any finite choice category defines a subspace in $U$ where all subsets are both choice categories and ``local'' sense categories.
\end{theorem}

In other words is the induced or relative topology on the subspace discrete both with respect to sense categories and choice categories. The proof is trivial\footnote{A Danish version of the proof can be found in \cite[1983, p. 372]{Mammen96}.}. The term ``correspondence'' refers to the fact that within any finite choice category is the logical structure reduced to the well-known classical ``Aristotelian'' logic in the same way as e.g. relativistic or quantum mechanical theories under limiting conditions are reduced to classical physics, which is Niels Bohr's famous correspondence principle, reciprocal to the abovementioned principle of conservative generalization.

\begin{theorem}(Th. 10, Globality): Any sense category in $U$ containing a non-empty choice category defines a subspace where all axioms Ax. 1-11 are satisfied.
\end{theorem}

The proof of this theorem is also trivial\footnote{A Danish version of the proof can be found in \cite[1983, p. 371]{Mammen96}.}.

The theorem tells that the structure defined by the 11 axioms is global or pervasive in $U$ and that it repeats itself in any detail as a fractal structure, or in mathematical terms that it is hereditary. It also says that the axiomatic system is rather ``immune'' to changes in definitions and interpretations of $U$ and its ``range''.

\subsection{The possible completeness of the axiomatic system}

The aim of the analysis until now has been to establish an understanding of the interface between man and the world of objects. This is of course not exhausting psychology in any way but just defining a foundation or basis on which to build an understanding of development of our cognition, actions, feelings, language, and much more.

In the present context we shall, however, dwell a little more on this basis itself, or in other words humans' immediate interface with the world of objects. One urgent problem is here if this basis, as defined by the 11 axioms, can be considered complete in the sense that there is not some third kind of category determining the structure of the interface. Could you from the 11 axioms deduce that such a third kind of category must necessarily exist?

Having Def. 1 in mind this question can be expressed in these two conjectures:

\medskip

\noindent {\bf Conjecture II.} (Claim of Completeness, CC).  
There exists a space on $U$ where any subset is a decidable category.
\medskip

or its negation:

\medskip
\noindent {\bf Conjecture III.}
There always exists a subset in $U$ which is not decidable.

\medskip

These opposing conjectures were put forward in \cite[1983, pp. 406-407]{Mammen96}. In \cite[1989, pp. xvi-xvii]{Mammen96} it was then proven that Conjecture III was true if the set of sense categories had a countable basis, but no more general proofs were established.

However, in 1994 J{\o}rgen Hoffmann-J{\o}rgensen proved in \cite{HJ2000} that if the sense categories had a basis with higher cardinality than $U$, then Conjecture II (CC) was true if Zorn's Lemma was true, Zorn's Lemma being equivalent with the Axiom of Choice (AC).

Hoffmann- J{\o}rgensen referred to the fact that Zorn's Lemma implied the existence of maximal perfect topologies \cite{Hewitt,Douwen} and proved that this existence further implied the Claim of Completeness (CC).

\medskip

Hoffmann-J{\o}rgensen then, as stated in the introduction, put forward the opposite implication as a conjecture:

\medskip
\noindent {\bf Conjecture I.} (Hoffmann-J{\o}rgensen; \cite[p. 86]{Mammen17}).
CC implies AC.
\medskip

It was rather surprising, that the set of 11 axioms combined with the claim of completeness seemed to imply an exotic structure as maximal perfect topologies, higher cardinalities, and perhaps also the axiom of choice. After all, taken separately the axioms were extremely simple as they were directly translatable into first-order-logic, and not more complicated than they could be explained on elementary school level.

However, Hoffmann-J{\o}rgensen, and in fact some colleagues in Moscow, were not able to prove the above conjecture.

\section{Mathematical background}

This section sets the general stage for the mathematical results of the paper, and collects various observations and lemmas about Mammen spaces (defined in the introduction), and the connection between Mammen spaces and maximal perfect Hausdorff topological spaces.

\subsection{General observations} Let $U$ be a non-empty set, let $(U,\mathcal S,\mathcal C)$ be a Mammen space with universe $U$, and let $\mathcal T$ be a perfect Hausdorff topology on $U$.

(1) If $I(\mathcal C)$ is the ideal generated by $\mathcal C$, then it is easy to verify that $(U,\mathcal S,I(\mathcal C))$ is also a Mammen space.

(2) If $I\neq\{\emptyset\}$ is an ideal on $U$ such that $I\cap\mathcal T=\{\emptyset\}$, then $(U,\mathcal T,I)$ is easily seen to be a Mammen space. 

(3) Since $\mathcal T$ is a perfect topology, (2) gives that $(U,\mathcal T,\fin(U))$ is a Mammen space, where $\fin(U)$ denotes the ideal of finite subsets of $U$.

(4) Let $\nd(\mathcal T)$ denote the ideal of nowhere dense sets in the topology $\mathcal T$. Then $(U,\mathcal T,\nd(\mathcal T))$ is a Mammen space. Note that $\nd(\mathcal T)\supseteq\fin$, so this gives us an example with a potentially richer family of choice categories.

\medskip

Recall from the introduction that a Mammen space $(U,\mathcal S,\mathcal C)$ is complete if every $X\subseteq U$ can be written as $X=S\cup C$, where $S\in\mathcal S$ and $C\in\mathcal C$. Building on (4) above, the next proposition tells us that in complete Mammen spaces, the sets $\nd(\mathcal T)$ are necessarily choice categories:

\begin{prop}\label{p.complete}
Suppose $(U,\mathcal S,\mathcal C)$ is a complete Mammen space. Then

(a) $X\subseteq U$ contains no non-empty open set in $\mathcal S$ iff $X\in\mathcal C$.

(b) $\mathcal C$ is an ideal, consisting precisely of the sets with empty interior.

(c) $\nd(\mathcal S)\subseteq\mathcal C$.
\end{prop}
\begin{proof}
(a) ``$\Longrightarrow$'': By completeness, $X=S\cup C$ for some $S\in\mathcal S$ and $C\in\mathcal C$. So if $X$ contains no non-empty set from $\mathcal S$, then $X=C$ follows.

``$\Longleftarrow$'': If $X\in\mathcal C$ and $S\subseteq X$ where $S\in\mathcal S$, then $X\cap S=S\in\mathcal C$ by (3.b) of Definition \ref{d.mammen}/Ax. 11, and $S=\emptyset$ follows by (3.a) of Definition \ref{d.mammen}/Ax. 6.

(b) $\mathcal C$ is closed under finite unions by (2.b)/Ax. 10. That $\mathcal C$ is closed under subsets is clear by (a), and (a) also gives that it consists precisely of the sets with empty interior.

(c) Clear by (b).
\end{proof}

The following simple combinatorial lemma will be used several times in the sequel; it was already observed by Jens Mammen in his early investigations of his axiom system, see e.g. \cite[1989 pp. xvi--xxi]{Mammen96}.

\begin{lemma}[Mammen]\label{l.lemmaM}
Let $(U,\mathcal S,\mathcal C)$ be a Mammen space. Suppose there is $X\subseteq U$ such that the following property holds:
\begin{quote}
$(*)$ For any $S\in\mathcal S\setminus\{\emptyset\}$ the sets $S\cap X$ and $S\setminus X$ are infinite.
\end{quote}
Then $(U,\mathcal S,\mathcal C)$ is not complete.
\end{lemma}
\begin{proof}
Suppose, seeking a contradiction, that $(U,\mathcal S,\mathcal C)$ were complete. Then
$$
X=S\cup C
$$
for some $S\in\mathcal S$ and $C\in\mathcal C$. By $(*)$, we can't have that $S\neq\emptyset$, so we must have $X=C$, so $X\in\mathcal C$. By the same reasoning, we must also have $U\setminus X\in\mathcal C$. But then $U=X\cup (U\setminus X)\in\mathcal C$ by (2.a) of Definition \ref{d.mammen}/Ax. 10; but this contradicts (3.a)/Ax. 6.
\end{proof}

From the previous Lemma, it is easy to derive the following:

\begin{theorem}[Mammen]
If $(U,\mathcal S,\mathcal C)$ is a complete Mammen space then $\mathcal S$ is not second countable.
\end{theorem}
\begin{proof}
If $(U_n)_{n\in\N}$ enumerates a countable basis for $\mathcal S$, then it is easy to construct from $(U_n)_{n\in\N}$ a set $X\subseteq U$ satisfying the property $(*)$.
\end{proof}

\begin{remark}
The reasoning of the previous proof can be adapted to prove a stronger result under the assumption of Martin's Axiom: Assuming $MA(\kappa)$, the basis of a $\mathcal S$ in a complete Mammen space must have cardinality $>\kappa$. See also the discussion of cardinal invariants and Theorem \ref{t.add} below.
\end{remark}

\subsection{Maximal perfect topologies and complete Mammen spaces} We now describe a method, due to Hoffmann-J{\o}rgensen, for obtaining complete Mammen spaces by considering \emph{maximal perfect topologies}.

\begin{definition}
Let $U$ be an infinite set. A perfect topology $\mathcal T$ on $U$ is said to be a  \emph{maximal perfect topology} if no topology finer than $\mathcal T$ is perfect.
\end{definition}

The next theorem is due to Hoffmann-J{\o}rgensen, \cite{HJ2000}. It provides the central connection between complete Mammen spaces and maximal perfect topologies. We note that the proof of this theorem (and the lemma following it) does not use the Axiom of Choice.

\begin{theorem}[Hoffmann-J{\o}rgensen]\label{t.hoffmann}
Let $\mathcal T$ be a perfect Hausdorff topology on an infinite set $U$. Then $\mathcal T$ is maximal if and only if every set $X\subseteq U$ can be written as $X=S\cup C$ where $S\in\mathcal T$ and $C$ is closed and discrete (and therefore is closed nowhere dense).
\end{theorem}

The following easy lemma will be used in the proof of Theorem \ref{t.hoffmann}, and in many other places later.

\begin{lemma}\label{l.maxtop}
Let $\calT$ be a perfect topology on a set $U\neq\emptyset$. Suppose $X\subseteq U$ is such that
\begin{equation}\label{eq.interscard}
(\forall V\in\calT)\  |X\cap V|\in\emptyinfinite,
\end{equation}
that is, $X\cap V$ is always either empty or infinite for every $V\in\mathcal T$. Then $\mathcal T\cup\{V\cap X:V\in\mathcal T\}$ is a basis of a perfect topology $\calT'\supseteq\calT$ with $X\in\calT'$.

It follows that $\calT$ is maximal if and only if every $X\subseteq U$ which satisfies Eq. (\ref{eq.interscard}) must be open. 

Moreover, if $\mathcal T$ is maximal, perfect and Hausdorff, then every discrete set is closed.
\end{lemma}

\begin{proof}[Proof of Lemma \ref{l.maxtop}]
It is clear that $\mathcal T\cup\{V\cap X:V\in\mathcal T\}$ is closed under finite intersections, and so forms a basis for a topology $\mathcal T'$ refining $\mathcal T$. That $\mathcal T'$ is perfect follows easily from Eq. (\ref{eq.interscard}).

For the ``moreover'', let $C\subseteq U$ be discrete; we will use the first part of the lemma to show that $U\setminus C$ is open. For this, suppose, seeking a contradiction, that $W\cap (U\setminus C)$ is finite and non-empty for some $W\in\mathcal T$. Since $\mathcal T$ is a perfect topology, every non-empty open set is infinite, and so $W\cap C$ must be infinite since $W\cap (U\setminus C)$ is finite. In particular, $W\cap C\neq\emptyset$; let $x\in W\cap C$. Since $C$ is discrete, we can find $V\in\mathcal T$ such that $V\cap C=\{x\}$. Now a contradiction with perfectness of $\mathcal T$ ensues, since $W\cap V=\{x\}$ is a finite non-empty open set.
 \end{proof}

\begin{proof}[Proof of Theorem \ref{t.hoffmann}]
``$\Longleftarrow$'': Suppose $X$ satisfies Eq. \ref{eq.interscard}. By Lemma \ref{l.maxtop}, we just need to prove that $X\in\mathcal T$. For this, write $X=S\cup C$ with $S\in\mathcal S$ and $C$ closed discrete. We may assume that $S\cap C=\emptyset$, since we can otherwise replace $S$ by the open set $S\cap(U\setminus C)$. We claim that $C=\emptyset$, and therefore $X=S\in\mathcal T$. Indeed, if $C\neq\emptyset$ were the case, let $x\in C$. Then, since $C$ is discrete, there would be $V\in\mathcal T$ such that $V\cap C=\{x\}$. Then, since $S\cap C=\emptyset$, we hace  $V\cap X=V\cap (S\cup C)=\{x\}$, which contradicts that $X$ satisfies Eq. (\ref{eq.interscard}).

``$\Longrightarrow$'': Let $X\subseteq U$, and let 
$$
C=\{x\in X: (\exists V\in\mathcal T)\ V\cap X=\{x\}\}.
$$
Clearly $C$ is discrete, and therefore closed by Lemma \ref{l.maxtop}. To see that $X\setminus C\in \mathcal T$, use Lemma \ref{l.maxtop}: If $(X\setminus C)\cap V$ was finite and non-empty for some $V\in\mathcal T$, then the Hausdorff property would give that $(X\setminus C)\cap V\subseteq C$.
\end{proof}

\begin{corollary}[Hoffmann-J{\o}rgensen \cite{HJ2000}]\label{c.hoffmann}
Let $\mathcal T$ be a maximal perfect topology on $U$, and let $\cd(\mathcal T)$ denote the family of closed discrete subset of $U$. Then

(a) $(U,\mathcal T,\cd(\mathcal T))$ is a complete Mammen space;

(b) $\nd(\mathcal T)=\cd(\mathcal T)=\{X\subseteq U: \interior(X)=\emptyset\}$.
\end{corollary}
\begin{proof}
(a) is clear by Theorem \ref{t.hoffmann}. (b) follows from Theorem \ref{t.hoffmann} and Proposition \ref{p.complete}(b).
\end{proof}

\subsection{Existence of maximal perfect topologies and complete Mammen spaces} A routine application of Zorn's lemma (and therefore AC) provides the following:

\begin{theorem}[Hewitt \cite{Hewitt}; uses Choice]\label{t.hewitt}
If $\mathcal T$ is a perfect topology on a set $U$, then there is a maximal perfect topology $\mathcal T'$ on $U$ such that $\mathcal T\subseteq\mathcal T'$.
\end{theorem}

Using this theorem and Corollary \ref{c.hoffmann}, we get:

\begin{corollary}[Hoffmann-J{\o}rgensen]
The Axiom of Choice implies that there are complete Mammen spaces. We can even obtain a complete Mammen space  with a countable universe.
\end{corollary}

\begin{proof}
The first statement is clear by Theorem \ref{t.hewitt} and Corollary \ref{c.hoffmann}. For the second part, take $U=\Q$ (the rationals), and extend the topology induced by open rational intervals to a maximal perfect topology.
\end{proof}





\section{Theorem B: Measurability and complete Mammen spaces}

In this section we will prove:

\begin{theorem}

(a) If all sets are Lebesgue measurable then there is no complete Mammen space with universe $\N$.

(b) If all sets are Baire measurable then there is no complete Mammen space with universe $\N$.
\end{theorem}

Of course, Solovay \cite{solovay70} famously showed that if ZF is consistent then so is ZF+``all sets are Lebesgue and Baire measurable''. (In the Lebesgue case we need an inaccessible cardinal to obtain this, but in the Baire case, Shelah \cite{shelah84} famously showed we don't) So the previous theorem tells us that at least \emph{some} amount of Choice is needed to obtain a complete Mammen space with a countably infinite universe. It is unclear if Lebesgue and Baire measurability has any influence on the existence of complete Mammen spaces with universes of higher cardinality than $\aleph_0$; see Question \ref{q.universelarge} later.

\begin{proof}
The proofs of (a) and (b) are virtually identical. We give the details for (a).

Identify $\mathcal P(\N)$ with $2^\N=\{0,1\}^\N$ in the natural way, and equip $2^\N$, and therefore $\mathcal P(\N)$, with the ``coin-flipping measure'' $\mu$, that is, the product measure on $\{0,1\}^\N$ where equal weight $1/2$ is given to $0$ and $1$. Then the function $\rho:\mathcal P(\N)\to\mathcal P(\N)$ defined by $\rho(A)=\N\setminus A$ is measure-preserving.

Let $(\N,\mathcal S,\mathcal C)$ be a Mammen space with universe $\N$, and assuming all subsets of $\mathcal P(\N)$ are $\mu$-measurable. We will show that $(\N,\mathcal S,\mathcal C)$ is not complete.

To see this, define for each $n\in\N$ the set
$$
\mathcal A_n=\{A\subseteq\N: (\forall V\in\mathcal S)\ n\in V\implies |V\cap A|=\infty\}.
$$
Note that $\mathcal A_n$ is ``$E_0$-invariant'', i.e., is invariant under finite changes: If $A\in\mathcal A_n$ and $B\subseteq\N$ is such that $A\triangle B$ is finite, then $B\in\mathcal A_n$. Since we're assuming that all sets are Lebesgue measurable, the $E_0$-invariance of $\mathcal A_n$ implies that $\mu(\mathcal A_n)=1$ or $\mu(\mathcal A_n)=0$.

\medskip

{\bf Claim:} $\mathcal P(\N)=\mathcal A_n\cup\rho(\mathcal A_n)$.
\medskip

{\it Proof of Claim}: Suppose not, and let $A\subseteq\N$ be such that $A\notin\mathcal A_n\cup\rho(\mathcal A_n)$. Then by definition of $\mathcal A_n$ there must be $V,V'\in\mathcal S$ such that $n\in V$ and $n\in V'$ and $A\cap V$ and $A^c\cap V'$ are finite. It follows that $V\cap V'\cap A$ and $V\cap V'\cap A^c$ are finite sets, and so $V\cap V'$ is finite. But since $n\in V\cap V'$, we have $V\cap V'\neq\emptyset$, contradicting that $\mathcal S$ is a perfect topology.\hfill\qed

\medskip

The previous claim gives that $\mu(\mathcal A_n)>0$ or $\mu(\rho(\mathcal A_n))>0$, but since $\mu(\rho(\mathcal A_n))=\mu(\mathcal A_n)$, it then follows $\mu(\mathcal A_n)>0$; and the $E_0$-invariance of $\mathcal A_n$ then gives us that $1=\mu(\mathcal A_n)=\mu(\rho(\mathcal A_n))$. It follows that
$$
\mu(\bigcap_{n\in\N} \mathcal A_n\cap \bigcap_{n\in\N} \rho(\mathcal A_n))=1,
$$
and so there is $X\in\bigcap_{n\in\N} \mathcal A_n\cap \bigcap_{n\in\N} \rho(\mathcal A_n)$. Then $X$ must be infinite, and $X\notin\mathcal S$.

Let $V\in\mathcal S\setminus\{\emptyset\}$, and let $n\in V$. Since $X\in \mathcal A_n\cap \rho(\mathcal A_n)$ we have $|V\cap X|=|V\cap X^c|=\infty$. So by Lemma \ref{l.lemmaM} $(\N,\mathcal S,\mathcal C)$ is not a complete Mammen space.
\end{proof}

The previous proof can easily be localized to pointclasses in Polish spaces (in the usual sense of descriptive set theory, see \cite{Kechris95}). In particular we have:

\begin{corollary}
There are no complete Mammen spaces $(U,\mathcal S,\mathcal C)$ where $\mathcal S$ is analytic as subsets of $\mathcal P(\N)$.
\end{corollary}
\begin{proof}
We just need to observe that if $\mathcal S$ is analytic, then $\mathcal A_n$ is co-analytic, and therefore Lebesgue measurable, and then the rest of the proof goes through unchanged.
\end{proof}

\begin{corollary}
If all sets are Lebesgue measurable (or all sets are Baire measurable) then there are no maximal perfect topologies on $\N$ (or any other countably infinite set).
\end{corollary}

\section{Theorem A: Maximal perfect topologies in $\hod^{V[G]}(A)$}

In this section, we will prove Theorem A. Specifically we will prove:

\begin{theorem}\label{t.mainA}
In the first Cohen model, every perfect topology can be extended to a maximal perfect topology.
\end{theorem}
Theorem A then follows by combining Theorem \ref{t.mainA} and Theorem \ref{t.hoffmann}.

\medskip

Our proof follows Repick\'y's \cite{Repicky} presentation of Halpern and L\'evy's theorem \cite{HaLe71} that the Boolean prime ideal theorem and the ultrafilter lemma (i.e., ``every ideal in a Boolean algebra can be extended to a prime ideal'' and ``every filter can be extended to an ultrafilter'', respectively) holds in the first Cohen model. In keeping with \cite{Repicky} and \cite{Jech03}, we will use $\omega$ for the set of non-negative integers.

\subsection{Notation and the first Cohen model} Our ground model will be called $V$. Let $\P\in V$ be the poset of all finite functions $p\subseteq (\omega\times\omega)\times\{0,1\}$. If $G\subseteq \P$ is a filter generic over $V$, then let $a:\omega\times\omega\to\{0,1\}$ be $a=\bigcup G$, and let $a_i(n)=a(i,n)$. Let $A=\{a_i:i\in\omega\}$. The ``first Cohen model'' is then $\hod^{V[G]}(A)$. It is well-known that $\hod^{V[G]}(A)$ is a model of ZF, but that the Axiom of Choice is false in $\hod^{V[G]}(A)$, since the set $A$ is infinite, yet $\hod^{V[G]}(A)$ believes that $A$ has no countable subsets. (An excellent, and brief, account of all the notions referred to in this paragraph can be found in \cite[Ch. 13-14]{Jech03}; a much fuller account of the first Cohen model can be found in \cite{Jech73}.)

Following Kechris \cite[Theorem 19.1]{Kechris95}, we will denote by $(C)^m$ the set of injective sequences in the set $C$ of length $m$. As always, $[C]^m$ denotes the set of all $m$-element subsets of $C$. We will use $C^m$ and $C^{<\omega}$ for the set of $m$-element sequences from $C$ and the set of all finite sequences (indexed from 0), respectively. (Repick\'y uses $^m C$ and $^{<\omega}C$ instead.)

\medskip

Next, we recall the two key lemmas from Repick\'y's paper:

\begin{lemma}[``Schema of continuity'', Lemma 2 in \cite{Repicky}]\label{l.rep1}\

Let $\varphi(w_1,\ldots, w_n, u,v)$ be a formula in the language of (ZF) set theory with free variables shown. Suppose that for some $x_1,\ldots, x_n\in V$, $m\in\omega$ and $s\in (A)^m$ we have $V[G]\models\varphi(x_1,\ldots, x_n,s,A)$. Then there is $k$ such that for any $\vec{t}\in A^m$ with $\vec{t}_i\supseteq s_i\restrict k$ for all $i<m$ we have
$$
V[G]\models \varphi(x_1,\ldots, x_n,\vec{t},A).
$$
Moreover, $k$ may be chosen such that the finite sequences $s_0\restrict k,\ldots,s_{m-1}\restrict k$ are pairwise incompatible.
\end{lemma}

Tracking Repick\'y's \cite{Repicky} again, we make the following definition:

\begin{definition}
Let $F\subseteq [A]^m$, and let $u_1,\ldots, u_m\in 2^{<\omega}$ be pairwise incompatible. We will say that $u_1,\ldots, u_m$ \emph{distinquish} $F$ if $F\cap N_{u_i}$ is a singleton for all $i\in\leq m$, where $N_{u_i}$ is the basic open neighbourhood in $2^{\omega}$ determined by $u_i$.
\end{definition}

\begin{lemma}[Corollary 3 in \cite{Repicky}]\label{l.rep2}
Let $\varphi(w_1,\ldots, w_n, v)$ be a formula in the language of set theory with free variables shown. Let $s\in A^{<\omega}$, $x_1,\ldots, x_n\in \od^{V[G]}[A,s]$, and let $F'\subseteq A\setminus \ran(s)$ be a finite set, and let $m=|F'|$.

Suppose $\varphi(x_1,\ldots, x_n, F')$ holds in ${V[G]}$. Then there are $u_1,\ldots, u_m\in 2^{<\omega}$ which distinguish the elements of $F'$, and $\varphi(x_1,\ldots, x_n,F)$ holds in $V[G]$ for any $F\in [A]^m$ such that $F$ is distinquished by $u_1,\ldots, u_m$.
\end{lemma}

\begin{proof}[Proof of Theorem \ref{t.mainA}] We will work in $V[G]$, so that $\od$ refers to $\od^{V[G]}$ and $\hod$ refers to $\hod^{V[G]}$, etc. Let $(X,\mathcal T)\in \hod(A)$ and suppose 
$$
\hod(A)\models\text{``$\mathcal T$ is a perfect topology on $X$''}.
$$
Then for some finite sequence $f\in A^{<\omega}$ we have $X,\mathcal T\in\od[A,f]$. For notational simplicity, we assume that $f=\emptyset$, that is, $X,\mathcal T\in\od[A]$, as the presence of the $f$ makes no difference for our argument.

There is a well-ordering of $\od[A]$ which itself is ordinal definable from $A$ (see \cite[Lemma 13.25]{Jech03}). Using this well-ordering, we can define a perfect topology $\mathcal T'\in\od[A]$ on $X$ with $\mathcal T'\supseteq \mathcal T$ which is maximal among perfect topologies ordinal definable from $A$. We claim that 
$$
\hod(A)\models\text{``$\mathcal T'$ is a maximal perfect topology''}.
$$
To see this, we will prove the following: 

\medskip

{\bf Claim:} If $\mathcal T'$ is maximal among perfect topologies on $X$ in $\od[A,s]$ for some $s\in A^{<\omega}$, then it is maximal among perfect topologies on $X$ in $\od[A,s\append a]$ for any $a\in A$. 

\medskip

If we can prove this claim, then an easy induction on $\lh(s)$ shows that $\mathcal T'$ is maximal among perfect topologies in $\od[A,s]$ for any $s\in A^{<\omega}$, and so $\mathcal T'$ is maximal in $\hod(A)$ (see e.g. \cite[pp. 186--188]{Jech03} for the general background on $\od$ and $\hod$).

\medskip

We now turn to the proof of the claim. As $s$ will play no role in our argument, we suppress it (that is, we give the argument for $s=\emptyset$, which is virtually identical to the argument for $s\neq\emptyset$). 

To prove the claim, we will use Lemma \ref{l.maxtop}. Suppose $a'\in A$ and $w\in\od[A,a']\cap \mathcal P(X)$, and for all $v\in\mathcal T'$ either $w\cap v=\emptyset$ or $w\cap v$ is infinite (in $V[G]$). By Lemma \ref{l.maxtop}, we need to show that $w\in\mathcal T'$. Assume for a contradiction that $w\notin\mathcal T'$. Let $\varphi$ be a formula such that
$$
w=\{x\in V[G]: V[G]\models\varphi(x,\alpha_1,\ldots,\alpha_n,a',A)\}.
$$
where $\alpha_1,\ldots,\alpha_n$ are ordinals. Then by Lemma \ref{l.rep2}, there is $u\in 2^{<\omega}$ such that for all $a\in A\cap N_u$ and all $v\in\mathcal T'$ we have either $w(a)\cap v=\emptyset$ or $w(a)\cap v$ is infinite, and $w(a)\notin\mathcal T'$, where
$$
w(a)=\{x\in V[G]: V[G]\models\varphi(x,\alpha_1,\ldots,\alpha_n,a,A\}.
$$
The set
$$
\mathcal S=\mathcal T'\cup\{v\cap \bigcap_{a\in F} w(a) : v\in\mathcal T'\wedge F\subseteq A \text{ finite}\wedge (\forall a\in F)\ u\subset a\}
$$
is definable from $\alpha_1,\ldots,\alpha_n$ and $A$, so is in $\od[A]$, and moreover it is the basis of a topology which is strictly finer than $\mathcal T'$. Since $\mathcal T'$ is maximal among perfect topologies in $\od[A]$, there must be some finite $F'\subseteq A$ and $v\in\mathcal T'$ such that $z(F')$ is finite, where in general we let
$$
z(F)=v\cap \bigcap_{a\in F} w(a).
$$
Note that by the assumptions on $u$, we must have $|F'|>1$. Let $m=|F'|$. We may assume that $m$ is minimal, that is, for no $E\subseteq A\cap N_u$ with $|E|<m$ do we have $z(E)$ finite.) 

Since $z(F')$ is a finite subset of $X\in \hod(A)$ we can find $s\in A^{<\omega}$ and $x_1,\ldots, x_n\in\od[A,s]$ and $V[G]\models z(F')=\{x_1,\ldots, x_n\}$. By Lemma \ref{l.rep2} we can find $u_1,\ldots, u_m\in 2^{<\omega}$ pairwise incompatible and extending $u$ which distinguish $F'$. Then for any $F\in [A]^m$ distinquished by $u_1,\ldots, u_m$ we have
$$
z(F)=\{x_1,\ldots, x_n\},
$$
which shows that $z(F')\in\od[A]$. 

Now, for each $1\leq i\leq m$, let
$$
y_i=\bigcup_{a\in A\cap N_{u_i}} w(a).
$$
Then $y_i\in \od[A]$, and $v\cap y_1\cap\cdots\cap y_m=\{x_1,\ldots, x_n\}$. We must have that the sets $v'\cap y_i$ are infinite or empty for any $v'\in\mathcal T'$ since this already holds for any $w(a)$ with $a\in A\cap N_u$. So $y_i\in\mathcal T'$ for all $i\leq m$ by the maximality of $\mathcal T'$ among perfect topologies in $\od[A]$. But now $v\cap y_1\cap\cdots\cap y_m$ is finite \emph{and} in $\mathcal T'$, which is impossible since $\mathcal T'$ is a perfect topology.
\end{proof}

\section{Some cardinal invariants}\label{s.cardinv}

In this section, and the next, we will study the cardinal invariants
\begin{align*}
\uT&=\inf\{\card(\mathcal T): \mathcal T \text{ is a maximal perfect Hausdorff topology on }\N\},\\
\uM&=\inf\{\card(\mathcal S): \mathcal S \text{ is a sense category of a complete Mammen space on }\N \}
\end{align*}
and we will prove items (1) and (2) of Theorem C. 

Since every maximal perfect topology gives rise to a Mammen space in a canonical way, we must have $\uM\leq \uT$. Of course, we must also have $\uM,\uT\leq 2^{\aleph_0}$ since any $\mathcal T$ and $\mathcal S$ above are subsets of $\mathcal P(\N)$. We do not know if $\uM<\uT$ is consistent with ZFC, see the question section at the end.

\medskip

Recall from \cite[p. 515]{Jech03} that $\add(\bp)$ denotes the \emph{additivity} of the ideal of meagre sets in Cantor space $2^\N$ (equivalently, in any Polish space), that is, $\add(\bp)$ is the least cardinal $\kappa$ such that the union of some family of $\kappa$ meagre sets is non-meagre.

\begin{theorem}\label{t.add}
If $(\N,\mathcal S,\mathcal C)$ is a complete Mammen space with universe $\N$, then $|\mathcal S|\geq\add(\bp)$.
\end{theorem}

\begin{proof}
Assume, aiming for a contradiction, that $|\mathcal S|<\add(\bp)$. Let 
$$
\mathcal S'=\{V\in\mathcal S: |V|=|\N\setminus V|=\aleph_0\}.
$$
This set is non-empty since $\mathcal S'$ is a perfect Hausdorff topology, and by assumption it has cardinality less than $\add(\bp)$. For each $V\in\mathcal S'$, let
$$
M_V=\{x\subseteq\N: |x\cap V|=|V\setminus x|=\aleph_0\}.
$$
Then $M_V$ is comeagre in $[\N]^\N$, and since $|\mathcal S'|<\add(\bp)$, the set $\bigcap_{V\in\mathcal S'} M_V$ is comeagre, and so non-empty. Let $x\in\bigcap_{V\in\mathcal S'} M_V$. Then no sense category is a subset of $x$, and so since $(\N,\mathcal S,\mathcal C)$ is a complete Mammen space, we must have $x\in\mathcal C$. Similarly, no sense category is a subset of $\N\setminus x$, so $\N\setminus x\in\mathcal C$. It follows that $\N=x\cup(\N\setminus x)\in\mathcal C$, contradicting that $\mathcal S\cap\mathcal C=\{\emptyset\}$.
\end{proof}

\begin{corollary}[Theorem C part (1)]
$\add(\bp)\leq\uM\leq \uT$.
\end{corollary}

\begin{corollary}[Theorem C part (2)]
Martin's Axiom (MA) implies that $\uM=\uT=2^{\aleph_0}$. So under MA,
the family of sense categories in a complete Mammen space always has cardinality $2^{\aleph_0}$. \end{corollary}
\begin{proof}
It is well-known (see \cite[Theorem 2.22]{kunen80}) that MA implies that $\add(\bp)=2^{\aleph_0}$. So by the previous corollary, MA implies that $2^{\aleph_0}\leq \uM\leq\uT$, and as noted above, $\uM\leq \uT\leq 2^{\aleph_0}$.
\end{proof}

\begin{remark} Since MA+$2^{\aleph_0}>\aleph_1$ is consistent (with ZFC, provided ZFC itself is consistent, see \cite[Theorem 6.3]{kunen80}) the previous corollary shows that it is consistent to have $2^{\aleph_0}>\aleph_1$ and no complete Mammen spaces on $\N$ with the set of sense categories having cardinality $\aleph_1$. In the next section, we show that it is also consistent to have $2^{\aleph_0}>\aleph_1$ \emph{and} have a complete Mammen space with the set of sense categories having cardinality $\aleph_1$.
\end{remark}

\section{$\uT$ and $\uM$ in the Baumgartner-Laver model}\label{s.bl}

By the \emph{Baumgartner-Laver model} $M[G]$ we mean the model of ZFC obtained by iteratively adding $\aleph_2$ Sacks reals to a model $M$, where $M$ satisfies the Continuum Hypothesis, CH. The purpose of this section is to prove

\begin{theorem}[Theorem C part (3)]\label{t.bltop}
In the Baumgartner-Laver model $M[G]$, there is a maximal perfect topology on $\N$ of cardinality $\aleph_1$. So in this model
$$
\aleph_1=\uT=\uM<2^{\aleph_0}=\aleph_2.
$$
\end{theorem}

In general, for a topology $\calT$ on $\N$ and $n\in\N$, let

\begin{equation}\label{eq.Tstar}
\calT^*(n)=\{V\setminus \{n\}: V\in\calT\wedge n\in V\}
\end{equation}

We will see that Theorem \ref{t.bltop} follows easily from Baumgartner and Laver's work once we prove:

\begin{theorem}\label{t.selectivetop}
Assume CH holds. Then there is a maximal perfect Hausdorff topology $\mathcal T$ on $\N$ such that for all $n\in\omega$, the family $\mathcal T^*(n)$ (defined in (\ref{eq.Tstar}) above) generates a selective ultrafilter on $\N\setminus\{n\}$.
\end{theorem}

\medskip

Before proving Theorem \ref{t.selectivetop}, we prove the following general lemma about perfect topologies in which $\calT^*(n)$ generates an ultrafilter.

\begin{lemma}\label{l.ultramax}
Let $\calT$ be a perfect topology on $\N$ such that for every $n\in\omega$, the set
$\calT^*(n)$ is the basis for an ultrafilter on $\N\setminus\{n\}$. Then $\calT$ is a maximal perfect topology.
\end{lemma}

\begin{proof}[Proof of Lemma \ref{l.ultramax}]
Let $X\subseteq\N$ and assume that
$$
(\forall V\in\calT)\  |X\cap V|\in\emptyinfinite.
$$
By Lemma \ref{l.maxtop} it is enough to show that $X\in\calT$. For this, it is enough to show that for any $n\in X$ there is $V\in\calT$ such that $n\in V\subseteq X$, since then
$$
X=\bigcup\{V\in\mathcal T: V\subseteq X\},
$$
which shows that $X\in\calT$.

So let $n\in X$. By the assumption on $\calT^*(n)$, there is $V\in\calT$ with $n\in V$ such that either $V\setminus\{n\}\subseteq X\setminus\{n\}$ or $(V\setminus\{n\})\cap (X\setminus\{n\})=\emptyset$. The latter can't be the case, since then $V\cap X=\{n\}$, which violates the assumption on $X$. So we must have $V\setminus\{n\}\subseteq X\setminus\{n\}$, from which it follows that $V\subseteq X$. 
\end{proof}

We now start working towards the main result of this section, Theorem \ref{t.bltop}. Recall the definition of a \emph{Ramsey ultrafilter} on $\N$ (see \cite[p. 71]{Jech03}):

\begin{definition}
A non-principal ultrafilter $\mathcal U$ on $\N$ is called a \emph{Ramsey} ultrafilter (also called a \emph{selective} ultrafilter) if for every partition $\{A_n:n\in\N\}$ of $\N$ into $\aleph_0$ pieces with each $A_n\notin\mathcal U$, there is $X\in\mathcal U$ such that $|A_n\cap X|\leq 1$ for all $n\in\N$.
\end{definition}

It is well-known, and quite easy, to show that if the Continuum Hypothesis (CH) holds, then there is a selective ultrafilter (which, since CH holds, must be of cardiality $\aleph_1$). Baumgartner and Laver, in their classic paper \cite{BL79}, showed the following:

\begin{theorem}[Baumgartner-Laver, 1979]\label{t.bl}
If $\aleph_2$ Sacks reals are added iteratively to a ground model $M$ which satisfies CH (the Continuum Hypothesis), then in the resulting model $M[G]$, it holds that $2^{\aleph_0}=\aleph_2$, and every selective ultrafilter in the ground model $M$ generates a selective ultrafilter in $M[G]$.

In particular, in $M[G]$ we have $\mathfrak u=\aleph_1<2^{\aleph_0}=\aleph_2$.
\end{theorem}

With Theorem \ref{t.bl} in mind, and Theorem \ref{t.bltop} in our sights, the next goal is to prove:

\begin{theorem}\label{t.selectivetop}
Assume CH holds. Then there is a maximal perfect Hausdorff topology $\mathcal T$ on $\N$ such that for all $n\in\omega$, the family $\mathcal T^*(n)$ (defined in (\ref{eq.Tstar}) above) generates a selective ultrafilter on $\N\setminus\{n\}$.
\end{theorem}

For the proof of Theorem \ref{t.selectivetop}, we need:

\begin{lemma}\label{l.addset}
Let $\mathcal T$ be a countable perfect Hausdorff topology on $\N$ and let $n\in\N$. Let $\mathcal A$ be a partition of $\N\setminus\{n\}$ into finitely or countably many pieces. Then there is an infinite set $B\subseteq\N\setminus\{n\}$ such that the following hold:
\begin{enumerate}
\item Either $B\subseteq A$ for some $A\in\mathcal A$, or $|B\cap A|\leq 1$ for all $A\in\mathcal A$.
\item $|B\cap V|\in\emptyinfinite$ for all $V\in\calT$.
\end{enumerate}

\end{lemma}
\begin{proof}[Proof of Lemma \ref{l.addset}]

The proof is divided into two cases.

\medskip

{\bf Case 1}: There is $\emptyset\neq V\in\calT$ such that $V\cap A\neq\emptyset$ for only finitely many $A\in\mathcal A$.

\smallskip

In this case there must be $A_1,\ldots, A_k\in\mathcal A$ such that
\begin{equation}\label{eq.cover}
V\setminus\{n\}\subseteq A_1\cup\ldots\cup A_k.
\end{equation}
Since there are only finitely many $A_1,\ldots,A_k$, it follows that there must be a non-empty $\tilde V\in\calT$ with $\tilde V\subseteq V$ such that for any non-empty $W\in\calT$ with $W\subseteq\tilde V$, we have
\begin{equation}\label{eq.partition}
\{i\in\{1,\ldots, k\}: |A_i\cap W|=\aleph_0\}=\{i\in\{1,\ldots, k\}: |A_i\cap \tilde V|=\aleph_0\}.
\end{equation}
Since $\tilde V$ is non-empty and $\calT$ is a perfect topology, $\tilde V$ must be infinite, and since $\tilde V\subseteq V$, it follows from (\ref{eq.cover}) that there is $i_0\leq k$ such that $|A_{i_0}\cap\tilde V|=\aleph_0$. Let $B=(A_{i_0}\cap \tilde V)\setminus\{n\}$. Then (1) in the lemma is clearly satisfied, and (2) is satisfied since for any $W\in\calT$ with $W\cap\tilde V\neq\emptyset$ we will have 
$$
|W\cap B|=|W\cap(A_{i_0}\cap \tilde V)\setminus\{n\}|=\aleph_0
$$
since the choice of $i_0$ and the fact that $\emptyset\neq W\cap \tilde V\in\mathcal T$ and $W\cap \tilde V\subseteq\tilde V$ ensures that $W\cap V\cap A_{i_0}$ is infinite.

\smallskip

{\bf Case 2}: For every $V\in\mathcal T\setminus\{\emptyset\}$ there are infinitely many $A\in\mathcal A$ such that $A\cap V\neq\emptyset$.

\smallskip

Let $E_{\mathcal A}$ denote the equivalence relation on $\N\setminus\{n\}$ corresponding to the partition $\mathcal A$, and let $[x]_{E_{\mathcal A}}$ denote the equivalence class of $x$. In the current case, each $V\in\calT\setminus\{\emptyset\}$ meets infinitely many $E_{\mathcal A}$-classes. Since there are only countably many $V\in\mathcal T$, an easy enumeration argument produces a family of sequences $(x_i^V)_{i\in\N, V\in\mathcal T\setminus\{\emptyset\}}$ such that
\begin{enumerate}
\item $x_i^V\in V\setminus\{n\}$ for all $i\in\N$;
\item The function $(i,V)\mapsto [x_i^V]_{E_\mathcal A}$ is injective from $\N\times(\mathcal T\setminus\{\emptyset\})$ into $(\N\cup\{n\})/E_{\mathcal A}$.
\end{enumerate}
Now let $B=\{x_i^V:(i,V)\in\N\times(\mathcal T\setminus\{\emptyset\})\}$. Then $|B\cap V|=\aleph_0$ for all $V\in\calT\setminus\{\emptyset\}$, and $|B\cap A|\leq 1$ for all $A\in\mathcal A$ by the injectivity of $(i,V)\mapsto [x_i^V]_{E_\mathcal A}$.
\end{proof}

\begin{proof}[Proof of Theorem \ref{t.selectivetop}]
Use CH to enumerate, for each $n\in\N$, all partitions (finite or infinite) of $\N\setminus\{n\}$ as $(\mathcal A_{n,\alpha})_{\alpha<\omega_1}$. Let $\mathcal T_0$ be a countable perfect Hausdorff topology on $\N$. For $n\in\N$ and $\alpha<\omega_1$, let
$$
\Lambda_{n,\alpha}=\{(i,\beta)\in\N\times\omega_1: (\beta=\alpha\wedge i<n)\vee(\beta<\alpha\wedge i\in\N)\}.
$$
By recursion on $\alpha<\omega_1$, we will define for each $n\in\N$ infinite sets $B_{n,\alpha}\subseteq \N\setminus\{n\}$ and perfect Hausdorff topologies $\calT_{n,\alpha}\supseteq\calT_0$ with the following properties:
\begin{enumerate}
\item $\calT_{n,\alpha}$ is the topology generated by
$$
\{B_{n,\alpha}\cup\{n\}\}\cup\calT_0\cup\bigcup_{(i,\beta)\in\Lambda_{n,\alpha}} \calT_{i,\beta},
$$
\item Either $B_{n,\alpha}\subseteq A$ for some $A\in\mathcal A_{n,\alpha}$, or $|B_{n,\alpha}\cap A|\leq 1$ for all $A\in\mathcal A_{n,\alpha}$.
\end{enumerate}
It is virtually clear by Lemma \ref{l.addset} that a recursion on $\alpha<\omega_1$ can be done: Having defined $\calT_{i,\beta}$ for all $(i,\beta)\in\Lambda_{n,\alpha}$, Lemma \ref{l.addset} can be applied with $\calT=\bigcup_{(i,\beta)\in\Lambda_{n,\alpha}}\calT_{i,\beta}$ to obtain $B_{n,\alpha}$ as desired, with (1) of Lemma \ref{l.addset} ensuring that $\calT_{n,\alpha}$ is a perfect topology, which is Hausdorff since $\calT_0\subseteq \calT_{n,\alpha}$.

\medskip

Let $\calT=\bigcup\{\calT_{n,\alpha}: n\in\N\wedge \alpha<\omega_1\}$. Then $\calT$ is a perfect Hausdorff topology in $\N$. To see that $\calT^*(n)$ generates a ultrafilter, let $A\subset\N\setminus\{n\}$, and let $\alpha$ be such that $\mathcal A_{n,\alpha}=\{A,A^c\}$. Then (2) guarantees that we must either have $B_{n,\alpha}\subseteq A$ or $B_{n,\alpha}\subseteq A^c$, while clearly $B_{n,\alpha}\in\calT^*(n)$ by (1). The selectivity property is also clear by (2). Finally, maximality of $\mathcal T$ follows from Lemma \ref{l.ultramax}. \end{proof}

Theorem \ref{t.bltop} is now an immediate corollary of Theorems \ref{t.selectivetop} and \ref{t.bl}.

\section{Open questions}

The following questions of a mathematical nature remain unsolved:

\subsection{Complete Mammen spaces and maximal perfect topologies} Inspired by Hoffmann-J\o rgensen, we have used maximal perfect topologies as a device to obtain complete Mammen spaces. It is natural to wonder how closely connected these two concepts are, specifically, we ask

\begin{question}
Does the existence of a complete Mammen space imply that there is a maximal perfect topology?
\end{question}

\subsection{First order compactness and Mammen spaces} One can quite easily make a first order formulation of Mammen's axiom system. The concept of completeness, though, is not so easily captured in such a first order axiomatization, since completeness of a space is a statement about \emph{all} subsets of the universe. Thus the following questions is natural:

\begin{question}\label{q.firstorder}
Does the first order compactness theorem imply that there is a complete Mammen space?\footnote{The second author of this paper at some point thought he had answered this question in the affirmative, and the first author announced this in reference \cite{Mammen19}. The second author dutifully retracts the claim of a solution, and the question remains wide open.} Does it imply there is a maximal perfect topology?
\end{question}

One may more generally ask:

\begin{question}\label{q.choice}
How weak a Choice principle is enough to ensure that a complete Mammen space exists?
\end{question}

Question \ref{q.firstorder} above can be thought of as a specific test case for the previous question.

\subsection{Regularity properties and the existence of complete Mammen spaces}

The next question takes aim at Question \ref{q.choice} from a different angle:

\begin{question}
Which regularity properties imply that there are no complete Mammen spaces with countable universe? E.g., if all sets are completely Ramsey, are there no complete Mammen spaces? What about Sacks, Miller, or Laver measurability, or other measurability notions that arise from arboreal forcing notions? (See e.g. \cite{BL99}.)
\end{question}

Of course, one may wonder if regularity properties have any influence on the existence of complete Mammen spaces with \emph{uncountable} universes; or if the existence of a complete Mammen space with uncountable universe can be achieved without appealing to Choice at all:

\begin{question}\label{q.universelarge}
Is it possible to prove in ZF without Choice (or with only weak Choice principes, such as countable choice or dependent choice) that there is a complete Mammen space with an uncountable universe?
\end{question}

\subsection{The cardinal invariants $\uM$ and $\uT$} We have seen in sections \ref{s.cardinv} and \ref{s.bl} the general inequalities
$$
\aleph_1\leq\add(\bp)\leq\uM\leq\uT\leq 2^{\aleph_0},
$$
and that (1) in models of Martin's Axiom, the last three $\leq$ are actually $=$, but (2) in the Baumgartner-Laver model, the first $\leq$ is actually $=$, and the last $\leq$ is actually $<$.

The most important unsolved question in this direction seems to be to separate $\uM$ and $\uT$:

\begin{question}
Is it consistent with ZFC to have $\mathfrak u_M<\mathfrak u_T$?
\end{question}

One may of course also wonder what the relation between $\mathfrak u_T$ and $\mathfrak u_M$ and the many other, well-known cardinal invariants that have been extensively studied. Most obviously, one may wonder what the connection between $\add(\lm)$, the additivity of the Lebesgue null ideal, and $\uM$ and $\uT$ is:

\begin{question}
Can $\add(\lm)\leq \uM$ be proven in ZFC?
\end{question}

Let us highlight one more question of this nature: Recall that $\mathfrak u$ denotes the smallest cardinality of a basis for a non-principal ultrafilter on $\N$. We ask:

\begin{question}
What is the relationship between $\mathfrak u$ and $\uM$ and $\uT$?
\end{question}

\section{Returning to Psychology}

The question of completeness of the basic interface, as described in the axioms 1-11, between human subjects and the world of objects is about the ultimate or ideal capacity of the interface. No human subject will be able to ``fill it out'' with categories realizing the complete case, and different people may differ in their ``repertoire'' of categories, and differ through their lives. The issue of completeness is therefore rather a question of whether sense and choice categories, or in short decidable categories, provide a sufficient conceptual frame, or system of reference, for describing people's factual system of categories and their development, e.g. in childhood, or if some third ``transcendental'' category should be needed by conceptual necessity, whether it is ``filled out'' or not.

The claim of completeness is therefore an expression of negation of \emph{a priori} limitations or restrictions on our access to subsets in the world of objects via decidable categories as defined by the axioms. But of course, there are also some factual limitations and restrictions of varying degree, which can be studied empirically within the complete frame. As human beings we are not reaching very far out in space and time, and there are many other practical restrictions on our activities.
However, there might still, hypothetically, be some further restrictions by principle on our factual ``repertoire'' of sense and choice categories.

\medskip

1. The completeness is dependent on some Choice principle weaker than the Axiom of Choice, but not yet made explicit, cf. Question 3 above.

It is therefore also an open question if this choice principle can be given an interpretation with some ``realism'', and accordingly some independent ``authority'', beyond the ``\emph{ad hoc}'' securing of completeness, or if it is too ``wild'' and should be replaced by a more ``modest'' choice principle, not securing completeness.
 
On the other hand, if completeness of the space is considered a sound and important principle in itself, we might have a criterion for deciding the corresponding choice principle to be fundamental, especially if it has further useful implications.

\medskip

2. Many psychological models of human perception and cognition, e.g. building on computer analogies and artificial intelligence, presuppose some degree of metric or regularity as basis for digital approximations or convergence towards our analogue reality. If such models are taken as premises it seems evident that completeness is excluded \emph{a priori}, cf. Question 4.

The same is the case if these models presuppose countable bases for their sense categories, which also excludes completeness, cf. Questions 5 and 6.

The last point raises the question of what is excluded \emph{a priori} working with computable or algorithmic models. It also raises the question of what is the reason for using algorithmic models of human activity at all. There is e.g. nothing in the function of the brain which points in that direction, despite popular ideas. It is true that some nerve impulses are of a binary on/off character, but they are occurring in continuous and not discrete time, and therefore not digital, but analogue, as the brain and body throughout.

However, already the fact that artificial intelligence models using pattern recognition are working exclusively on sense categories, although often within a user-defined finite frame of names referring to choice categories (e.g. persons or places), means that they \emph{a priori} are non-complete. Further, like a book, the AI models don't know the referents of the names. That is the user's human privilege.

This does not mean that algorithmic models can't be used as tools modeling domains with some regularity properties and being digitalized by intelligent humans. But they can't model humans themselves and their relations to the world, not even approximately.

\bibliographystyle{alpha}
\bibliography{stmodelmind}

\end{document}